\documentclass[a4paper,11pt,oneside]{amsart}

\usepackage{amsmath,amssymb}
\usepackage{microtype}
\usepackage[a4paper,left=25mm,right=25mm,top=26mm,textheight=253mm,headheight=0pt,headsep=0pt,footskip=10mm]{geometry}
\usepackage{xcolor}
\definecolor{citationred}{RGB}{128,0,0}
\definecolor{referenceblue}{RGB}{0,0,128}
\usepackage[colorlinks=true,citecolor=citationred,linkcolor=referenceblue,urlcolor=black,bookmarksnumbered=true]{hyperref}
\usepackage{bookmark}
\usepackage{etoolbox}
\usepackage{enumitem}
\usepackage{needspace}
\setlist[enumerate,1]{label=\textup{(\roman*)},ref=\roman*,leftmargin=2.5em,topsep=2pt,partopsep=0pt,noitemsep}

\AtBeginDocument{\setlength{\parskip}{0.9pt plus .3pt minus .3pt}}

\AtBeginEnvironment{thebibliography}{\setlength{\parskip}{0pt}\raggedright}

\allowdisplaybreaks[1]
\makeatletter
\@namedef{subjclassname@2020}{\textup{2020} Mathematics Subject Classification}
\patchcmd{\@maketitle}{\global\topskip42\p@}{\global\topskip20\p@}{}{\PackageError{manuscript}{Cannot adjust title position}{}}
\patchcmd{\@maketitle}{\dimen@34\p@}{\dimen@22\p@}{}{\PackageError{manuscript}{Cannot adjust title spacing}{}}
\patchcmd{\abstract}{\leftmargin3pc}{\leftmargin2pc}{}{\PackageError{manuscript}{Cannot adjust abstract width}{}}
\makeatother

\newtheorem{theorem}{Theorem}[section]
\newtheorem{maintheorem}{Theorem}

\newtheorem{proposition}[theorem]{Proposition}
\newtheorem{corollary}[theorem]{Corollary}
\newtheorem{conjecture}[theorem]{Conjecture}
\theoremstyle{definition}
\newtheorem{question}[theorem]{Question}
\newtheorem{example}[theorem]{Example}
\theoremstyle{remark}

\theoremstyle{plain}

\DeclareMathOperator{\Hom}{Hom}
\DeclareMathOperator{\Ext}{Ext}
\DeclareMathOperator{\tr}{tr}
\newcommand{\kk}{\Bbbk}
\newcommand{\m}{\mathfrak m}
\newcommand{\p}{\mathfrak p}

\title[Nearly Gorenstein rings and products]{Ideals Defining Nearly Gorenstein Rings\\Are (Still Almost) Never Products}
\author{Sora Miyashita}
\address{Department of Pure and Applied Mathematics, Graduate School of Information Science and Technology, The University of Osaka, Suita, Osaka 565-0871, Japan}
\email{u804642k@ecs.osaka-u.ac.jp}
\date{October 2026}
\hypersetup{pdftitle={Ideals Defining Nearly Gorenstein Rings Are (Still Almost) Never Products},pdfauthor={Sora Miyashita},pdfsubject={Commutative algebra}}

\begin{document}

\begin{abstract}
We prove that a graded ideal of height at least two in a positively graded polynomial ring over a field cannot be a product of two proper graded ideals if it defines a Cohen--Macaulay, generically Gorenstein ring. For a nonzero proper graded ideal, the ring defined by its square or any higher power has a nilpotent canonical trace ideal if and only if the ideal has height at least two. In characteristic zero, when the ideal has height at least two and the resulting ring is Cohen--Macaulay, its canonical trace ideal lies in the power of its graded maximal ideal whose exponent is one less than the defining exponent. Equality holds when the ideal is the graded maximal ideal of the polynomial ring.

\end{abstract}

\maketitle
\thispagestyle{plain}

\setlength{\parskip}{4pt plus .3pt minus .3pt}
\section{Introduction}\label{sec:introduction}

Throughout, $S=\kk[x_1,\ldots,x_n]$ is a positively graded polynomial ring over a field, with $\deg x_i>0$ and $\m=(x_1,\ldots,x_n)$. For a positively graded $\kk$-algebra $R$ with $R_0=\kk$, let $\m_R$ denote its graded maximal ideal. A Noetherian ring $R$, not necessarily Cohen--Macaulay, is \emph{generically Gorenstein} if $R_{\p}$ is Gorenstein for every minimal prime $\p$ of $R$.

Huneke~\cite[Theorem~2]{Huneke} proved that an ideal of height at least two in an unramified regular local ring cannot be a product of two proper ideals if it defines a Gorenstein ring. We prove the following extension for graded ideals in $S$.

\begin{maintheorem}\label{thm:cm-products}
Let $I\subsetneq S$ be a graded ideal of height at least two. If $S/I$ is Cohen--Macaulay and generically Gorenstein, then $I$ cannot be a product of two proper graded ideals.
\end{maintheorem}

We call $R$ \emph{nearly Gorenstein} if $R$ is Cohen--Macaulay and $\tr_R(\omega_R)\supseteq\m_R$ (see \cite[Definition~2.2]{HHS}). If $\dim R>0$, every minimal prime of $R$ is strictly contained in $\m_R$, since the minimal primes are homogeneous. At every such prime, the canonical trace ideal of a nearly Gorenstein ring localizes to the unit ideal, so the localization is Gorenstein. In particular, if $S/IJ$ is nearly Gorenstein of positive dimension and $I,J$ are nonzero proper graded ideals, then each of $I$ and $J$ is generated by one element (Corollary~\ref{cor:positive-nearly}). The same conclusion holds if $S/IJ$ is almost Gorenstein graded (Corollary~\ref{cor:almost-products}).

In dimension zero, even two distinct primary ideals can have a product defining a nearly Gorenstein ring (Example~\ref{ex:products}). When $n=2$, the ring $S/\m^k$ is nearly Gorenstein if and only if $k\leq2$ (see \cite[proof of Theorem~3.12]{GHHM}). We ask whether the same bound on the exponent holds for every graded ideal of height at least two.

\begin{question}\label{question:powers}
Let $I\subsetneq S$ be a graded ideal of height at least two. Can $S/I^k$ be nearly Gorenstein for any $k\geq3$?
\end{question}

Our second main result answers Question~\ref{question:powers} in characteristic zero by bounding the canonical trace ideal.

\begin{maintheorem}\label{thm:powers}
Suppose that $\operatorname{char}\kk=0$. Let $I\subsetneq S$ be a graded ideal of height at least two and $k\geq2$. If $R=S/I^k$ is Cohen--Macaulay, then $\tr_R(\omega_R)\subseteq\m_R^{k-1}$. Equality holds for $I=\m$. In particular, $R$ is not nearly Gorenstein when $k\geq3$.
\end{maintheorem}

In positive dimension the answer is negative over any field (Corollary~\ref{cor:trace-powers}). Conjecture~\ref{conj:positive-characteristic} concerns the Artinian case in positive characteristic.

\setlength{\parskip}{0.9pt plus .3pt minus .3pt}
\section{Products and canonical trace ideals}\label{sec:products}

In this section we prove Theorem~\ref{thm:cm-products}. We first recall the definition of the canonical module and its trace ideal.

Let $R$ be a Noetherian positively graded $\kk$-algebra with $R_0=\kk$. Its graded canonical module is $\omega_R=H^{\dim R}_{\m_R}(R)^\vee$ (see \cite[Definition~(2.1.2)]{GW}), where $M^\vee=\bigoplus_{j\in\mathbb Z}\Hom_{\kk}(M_{-j},\kk)$ is the graded $\kk$-dual. For $R=S/J$, one has $\omega_R\cong\Ext_S^{\operatorname{ht}J}(R,\omega_S)$, where $\omega_S=S(-\sum_i\deg x_i)$ (see \cite[Corollary~(2.2.8) and Proposition~(2.2.9)]{GW}).

The canonical trace ideal is $\tr_R(\omega_R)=\sum_{\varphi\in\Hom_R(\omega_R,R)}\varphi(\omega_R)$, where the sum includes homogeneous maps of every degree. If $R$ is Cohen--Macaulay, then $\tr_R(\omega_R)=R$ if and only if $R$ is Gorenstein (see \cite{HK} and \cite[Lemma~2.1]{HHS}).

The canonical trace ideal contains a nonzerodivisor if and only if $R$ is unmixed and generically Gorenstein (see \cite[Remarks~2.9~(3) and~2.13]{KM}). Here \emph{unmixed} means that $\dim R/\p=\dim R$ for every $\p\in\operatorname{Ass}R$. If $\dim R>0$ and $\tr_R(\omega_R)\supseteq\m_R$, then $R$ is generically Gorenstein (see \cite[Remark~2.9~(6)]{KM}).

We prove a stronger form of Theorem~\ref{thm:cm-products}, with $(S_2)$ in place of Cohen--Macaulayness.

\begin{theorem}\label{thm:s2-products}
Let $I,J\subsetneq S$ be nonzero graded ideals such that $R=S/IJ$ satisfies $(S_2)$.
\begin{enumerate}
\item\label{cond:s2-products-primes} If either factor has height at least two, every prime $\p$ minimal over $I+J$ is minimal over $IJ$, and $R_{\p}$ is not Gorenstein.
\item\label{cond:s2-products-principal} The ring $R$ is generically Gorenstein if and only if both $I$ and $J$ are principal.
\end{enumerate}
\end{theorem}

\begin{proof}[Proof of Theorems~\ref{thm:cm-products} and~\ref{thm:s2-products}]
We first prove \textup{(\ref{cond:s2-products-primes})}. Assume $\operatorname{ht}I\geq2$ or $\operatorname{ht}J\geq2$ and let $\p$ be minimal over $I+J$. Suppose $\p$ is not minimal over $IJ$. Localize at $\p$ and retain the notation. Then $I+J$ is primary to the maximal ideal $\m$ of $S$ and $\dim R>0$. Since $(I+J)(I\cap J)\subseteq IJ$, a power of $\m$ annihilates $(I\cap J)/IJ$. This module lies in $H^0_{\m}(R)$, which vanishes by $(S_2)$. Thus $I\cap J=IJ$ and $\operatorname{Tor}_1^S(S/I,S/J)=0$ (see \cite[Exercise~A3.17]{Eisenbud}).

Since $S$ is regular local, the vanishing of $\operatorname{Tor}_1$ gives $\operatorname{Tor}_i^S(S/I,S/J)=0$ for all $i>0$ by \cite[Corollary~1]{Lichtenbaum}. The depth formula (see \cite[Theorem~1.2]{Auslander}) gives
\[
\operatorname{depth}(S/I)+\operatorname{depth}(S/J)
=\operatorname{depth}S+\operatorname{depth}\bigl((S/I)\otimes_S(S/J)\bigr)=\dim S,
\]
where the last equality holds because $S$ is regular and $(S/I)\otimes_S(S/J)\cong S/(I+J)$ is Artinian. Since $S$ is a domain and $I,J\ne0$, both $S/I$ and $S/J$ have dimension strictly less than $\dim S$. Their depths are therefore strictly less than $\dim S$, and the equality forces both depths to be positive.

As $I\cap J=IJ$, there is an exact sequence $0\to R\to S/I\oplus S/J\to S/(I+J)\to0$. The associated local cohomology sequence (see \cite[Section~3.5]{BH}) contains
\[
0=H^0_{\m}(S/I\oplus S/J)\longrightarrow H^0_{\m}(S/(I+J))\longrightarrow H^1_{\m}(R).
\]
Since $S/(I+J)$ has finite length, the middle term is $S/(I+J)$ itself. It injects into $H^1_{\m}(R)$, so $H^1_{\m}(R)\ne0$ and $(S_2)$ gives $\dim R=1$. The inclusions $IJ\subseteq I,J$ give $\dim(S/I),\dim(S/J)\leq\dim R=1$. Both rings have positive depth, so each has depth and dimension one. The depth formula now gives $\dim S=2$. In this two-dimensional local domain, $I$ and $J$ have height one because they are nonzero and $\dim(S/I)=\dim(S/J)=1$.

Return to the original $S$. A height-one prime over $IS_{\p}$ contracts to a height-one prime over $I$ because localization preserves all chains below a prime contained in $\p$. The same argument applies to $J$. This contradicts $\operatorname{ht}I\geq2$ or $\operatorname{ht}J\geq2$.

Now $\p$ is minimal over $IJ$. Hence $\sqrt{(IJ)S_{\p}}=\p S_{\p}$ and $\operatorname{ht}_{S_{\p}}((IJ)S_{\p})=\dim S_{\p}=\operatorname{ht}_S\p\geq2$. For the inequality, $\p$ contains both $I$ and $J$. Also $IS_{\p}$ and $JS_{\p}$ are proper. The regular local ring $S_{\p}$ contains $\kk$ and is unramified. Thus $S_{\p}/((IS_{\p})(JS_{\p}))\cong R_{\p}$ is not Gorenstein by \cite[Theorem~2]{Huneke}.

For \textup{(\ref{cond:s2-products-principal})}, suppose $R$ is generically Gorenstein. The sum $I+J\subseteq\m$ has a minimal prime, so \textup{(\ref{cond:s2-products-primes})} gives $\operatorname{ht}I=\operatorname{ht}J=1$. Every height-one prime of $S$ is generated by an irreducible polynomial. Since $I$ and $J$ are each contained in such a prime, the greatest common divisors $f$ and $g$ of their homogeneous generators are nonunits. Let $I'$ be generated by the quotients of the generators of $I$ by $f$. Divide the generators of $J$ by $g$ in the same way to obtain $J'$. Then $I=fI'$ and $J=gJ'$. The generators of $I'$ have greatest common divisor $1$, so no height-one prime contains $I'$. The same holds for $J'$.

If $I'J'\ne S$, then $\operatorname{ht}(I'J')\geq2$ since any prime containing $I'J'$ contains $I'$ or $J'$. The product $IJ=(fg)(I'J')$ defines the same $(S_2)$ ring $R$. Both factors are proper and their sum is contained in $\m$. Choose $\p$ minimal over $(fg)+I'J'$. By \textup{(\ref{cond:s2-products-primes})}, $\p$ is minimal over $IJ$ and $R_{\p}$ is not Gorenstein. This contradicts generic Gorensteinness, so $I'J'=S$. The inclusion $I'J'\subseteq I'\cap J'$ gives $I'=J'=S$. The converse holds because $S/IJ$ is a hypersurface when both factors are principal.

Theorem~\ref{thm:cm-products} follows from \textup{(\ref{cond:s2-products-principal})} because a Cohen--Macaulay ring satisfies $(S_2)$ and a nonzero principal ideal of $S$ has height one.
\end{proof}

We first apply the theorem to nearly Gorenstein rings of positive dimension.

\begin{corollary}\label{cor:positive-nearly}
Let $I,J\subsetneq S$ be nonzero graded ideals. If $S/IJ$ is nearly Gorenstein and has positive dimension, then $I$ and $J$ are each generated by one element.
\end{corollary}

\begin{proof}
Apply Theorem~\ref{thm:s2-products}\textup{(\ref{cond:s2-products-principal})}, since $S/IJ$ satisfies $(S_2)$ and is generically Gorenstein.
\end{proof}

A Cohen--Macaulay graded ring $R$ is \emph{almost Gorenstein} if it admits a graded exact sequence $0\to R\to\omega_R(-a(R))\to C\to0$ in which $C=0$ or the minimal number of generators of $C$ equals its Hilbert--Samuel multiplicity with respect to $\m_R$. Here $a(R)=-\min\{j:[\omega_R]_j\ne0\}$ (see \cite[Definition~8.1]{GTT}).

Theorem~\ref{thm:s2-products} gives the following analogue of Corollary~\ref{cor:positive-nearly} for almost Gorenstein graded rings.

\begin{corollary}\label{cor:almost-products}
Let $I,J\subsetneq S$ be nonzero graded ideals. If $S/IJ$ is an almost Gorenstein graded ring, then $I$ and $J$ are each generated by one element. In particular, a one-dimensional almost Gorenstein graded ring $S/I^k$ with $k\geq2$ is a hypersurface.
\end{corollary}

\begin{proof}
Let $R=S/IJ$ and let $\p$ be a minimal prime of $R$. Localizing the embedding $R\to\omega_R(-a(R))$ gives $R_{\p}\hookrightarrow(\omega_R)_{\p}$. The ring $R_{\p}$ is Artinian and its canonical module has length $\ell(R_{\p})$. Thus the embedding is an isomorphism and $R_{\p}$ is Gorenstein (see \cite[Lemma~3.1]{GTT}). Apply Theorem~\ref{thm:s2-products}\textup{(\ref{cond:s2-products-principal})} to the Cohen--Macaulay ring $R$. For the last assertion set $J=I^{k-1}$.
\end{proof}

The condition of positive dimension in Corollary~\ref{cor:positive-nearly} cannot be omitted.

\begin{example}\label{ex:products}
Let $S=\kk[x,y]$, $I=(x,y^a)$ and $J=(x^b,y)$ with $a,b\geq2$. Then $IJ=(xy,x^{b+1},y^{a+1})$ and $R=S/IJ$ is the fiber product of $\kk[x]/(x^{b+1})$ and $\kk[y]/(y^{a+1})$ over $\kk$. In this case $\operatorname{Soc}R=\kk x^b\oplus\kk y^a$. To compute the trace, dualize the inclusion $xR\hookrightarrow R$ over $\kk$. This gives a surjection $\omega_R\to(xR)^\vee$. Up to shift, $xR\cong R/(y,x^b)\cong\kk[x]/(x^b)$ is isomorphic to its $\kk$-dual. Thus $xR$ is the image of a map $\omega_R\to R$, and the same holds for $yR$. This gives $\m_R=xR+yR\subseteq\tr_R(\omega_R)$. The socle has dimension two, so $R$ is not Gorenstein and $\tr_R(\omega_R)=\m_R$ (see also \cite[Theorem~A~(i)~(a)]{Miyashita}).
\end{example}

\Needspace{4\baselineskip}
\section{Powers and canonical trace ideals}\label{sec:powers}

In this section we prove Theorem~\ref{thm:powers} and compute the canonical trace ideal in height one.

\begin{proof}[Proof of Theorem~\ref{thm:powers}]
Let $u\in S$ represent a nonzero element of $\tr_R(\omega_R)$. Let $a=\max\{r\geq0:I\subseteq\m^r\}$ and choose $g\in I\setminus\m^{a+1}$. Denote the homogeneous parts of least ordinary total degree of $u$ and $g$ by $u_d$ and $g_a$. Setting $x_n=1$ is injective on ordinary homogeneous polynomials of a fixed degree, since that degree and the exponents of $x_1,\ldots,x_{n-1}$ determine the exponent of $x_n$. The field $\kk$ is infinite, so there are $c_1,\ldots,c_{n-1}\in\kk$ at which both $u_d(x_1,\ldots,x_{n-1},1)$ and $g_a(x_1,\ldots,x_{n-1},1)$ are nonzero.

Localize at $\m$ and retain the notation $S,I,R$. The canonical module localizes, and its trace ideal localizes by \cite[Proposition~\mbox{2.8(viii)}]{Lindo}. In particular the class of $u$ still belongs to the trace. Set $P=(x_i-c_i x_n\mid i<n)$. Before localization, imposing $x_i=c_i x_n$ for $i<n$ leaves $\kk[x_n]$ and sends $\m$ to $(x_n)$. Therefore $D:=S/P\cong\kk[z]_{(z)}$ with $x_n$ mapping to $z$ and $x_i$ to $c_i z$ for $i<n$.

For $v$ in the original polynomial ring, let $v(z)$ denote its image in $D$. Every term of ordinary total degree $r$ maps to a scalar multiple of $z^r$. Our choice of the $c_i$ gives $u(z)\in z^dD\setminus z^{d+1}D$ and $g(z)\in z^aD\setminus z^{a+1}D$. Every element of the original $I$ maps into $z^aD$ and $g(z)$ generates $z^aD$. Thus $ID=z^aD$ and $R/PR\cong D/(z^{ak})$.

Suppose $z^{ak-k+1}(\omega_R/P\omega_R)=0$. Every $R$-linear map $\omega_R\to R$ sends $P\omega_R$ into $PR$, so the induced map $\omega_R/P\omega_R\to R/PR$ has image annihilated by $z^{ak-k+1}$. The residue of $u(z)$ is in the sum of these images since the class of $u$ is in the trace. It follows that $z^{ak-k+1}u(z)\in z^{ak}D$, so $u(z)\in z^{k-1}D$. Thus $d\geq k-1$, as required. It remains to prove the asserted annihilation.

Let $h=\operatorname{ht}I$ and set $M_j=\operatorname{Tor}_{h-1}^S(S/I^j,D)$. Choose $f\in I$ with image $z^a$ in $D$. Multiplication by $f^j$ on a free resolution of $S/I^j$ lifts the zero map on $S/I^j$. After tensoring the resolution with $D$, this chain map is multiplication by $z^{aj}$. Thus $z^{aj}$ annihilates $\operatorname{Tor}_i^S(S/I^j,D)$ for every $i$.

Since $R$ is Cohen--Macaulay over the regular local ring $S$, take a finite free $S$-resolution $F_\bullet$ of $R$ of length $h$ and set $C_\bullet=F_\bullet\otimes_S D$. We have $H_h(C_\bullet)\subseteq C_h$ and $z^{ak}H_h(C_\bullet)=0$, so $H_h(C_\bullet)=0$ since $C_h$ is free over $D$. On the other hand $\omega_R=\Ext^h_S(R,S)$ is the cokernel at the end of the $S$-dual of $F_\bullet$. Reducing this presentation modulo $P$ gives $\omega_R/P\omega_R\cong H^h(\Hom_D(C_\bullet,D))$ since $F_\bullet$ is finite free. The universal coefficient theorem over $D$ gives (see \cite[Theorem~3.6.5]{Weibel})
\[
0\longrightarrow\Ext^1_D(M_k,D)\longrightarrow\omega_R/P\omega_R\longrightarrow\Hom_D(H_h(C_\bullet),D)\longrightarrow0.
\]

Thus $\omega_R/P\omega_R\cong\Ext^1_D(M_k,D)$. Over the discrete valuation ring $D$, the finite torsion module $M_k$ is a direct sum of modules $D/(z^b)$. The resolution $0\to D\xrightarrow{z^b}D\to D/(z^b)\to0$ gives $\Ext^1_D(D/(z^b),D)\cong D/(z^b)$. Hence $\omega_R/P\omega_R\cong M_k$ as $D$-modules.

For $j\geq2$ the quotient $S/I^j\to S/I^{j-1}$ gives the upper arrow in the diagram below. The vertical arrows come from $D\to D/zD=\kk$.
\[
\begin{array}{ccc}
M_j & \xrightarrow{\alpha_j} & M_{j-1}\\
\big\downarrow{\scriptstyle\ \rho_j} && \big\downarrow{\scriptstyle\ \rho_{j-1}}\\
\operatorname{Tor}_{h-1}^S(S/I^j,\kk) & \xrightarrow{\quad 0\quad} & \operatorname{Tor}_{h-1}^S(S/I^{j-1},\kk).
\end{array}
\]
To check the lower arrow, work in the original polynomial ring $\kk[x_1,\ldots,x_n]$. For $f_1,\ldots,f_j$ in the original ideal $I$, the product rule gives
\[
\frac{\partial(f_1\cdots f_j)}{\partial x_i}=\sum_{r=1}^j\frac{\partial f_r}{\partial x_i}\prod_{s\ne r}f_s\in I^{j-1}.
\]
Thus $\partial(I^j)\subseteq I^{j-1}$, where $\partial(I^j)$ is generated by the partial derivatives of elements of $I^j$. Since $I^j\subseteq (x_1,\ldots,x_n)^2$, \cite[Proposition~3.5]{Maleki} makes the lower map zero in degree $h-1>0$. This remains true after localization. The long exact Tor sequence for $0\to D\xrightarrow{z}D\to\kk\to0$ gives $\ker\rho_{j-1}=zM_{j-1}$. The square gives $\rho_{j-1}\alpha_j=0$, so $\alpha_j(M_j)\subseteq zM_{j-1}$.

Multiplication by $f$ gives $S/I^{j-1}\to S/I^j$ because $fI^{j-1}\subseteq I^j$. Let $\beta_j\colon M_{j-1}\to M_j$ be the corresponding map. Composing $S/I^j\to S/I^{j-1}$ with multiplication by $f$ gives multiplication by $f$ on $S/I^j$. Thus $\beta_j\alpha_j=z^a\operatorname{id}_{M_j}$.

For $w\in M_j$ take $w'\in M_{j-1}$ such that $\alpha_j(w)=zw'$. Then $z^aw=z\beta_j(w')$. In particular, if $z^sM_{j-1}=0$ with $s\geq1$, then $z^{s+a-1}M_j=0$. Starting with $z^aM_1=0$, we get $z^{a+(j-1)(a-1)}M_j=0$ for each $j\geq1$. Taking $j=k$ gives $z^{ak-k+1}M_k=0$. The isomorphism with $M_k$ now gives $z^{ak-k+1}(\omega_R/P\omega_R)=0$.

The claimed containment follows. If $I=\m$, then $\operatorname{Soc}R=\m_R^{k-1}$. For any $s\in\operatorname{Soc}R$ compose a surjection $\omega_R\to R/\m_R$ with the map $R/\m_R\to R$ sending $1$ to $s$. This puts $s$ in $\tr_R(\omega_R)$ and proves equality (see also \cite[Lemma~3.1]{Kumashiro}).
\end{proof}

In height one the canonical trace ideal has an explicit description in any characteristic.

\begin{proposition}\label{prop:height-one}
Let $I\subsetneq S$ be a nonzero graded ideal and $k\geq1$, and set $R=S/I^k$. If $\operatorname{ht}I=1$ and $f$ is a greatest common divisor of the generators of $I$, then $\tr_R(\omega_R)=(I:_S f)^k/I^k$. For $k\geq2$, the following conditions are equivalent:
\begin{enumerate}
\item $\operatorname{ht}I=1$;
\item $R_{\p}$ is Gorenstein for some minimal prime $\p$ of $R$;
\item the canonical trace ideal of $R$ is not nilpotent.
\end{enumerate}
\end{proposition}

\begin{proof}
For $k\geq2$ we first show \textup{(iii)} $\Rightarrow$ \textup{(ii)} $\Rightarrow$ \textup{(i)}. If the trace is not nilpotent, then it is not contained in some minimal prime $\p$ of $R$. It follows that $\tr_{R_{\p}}((\omega_R)_{\p})=R_{\p}$. By \cite[Remark~2.9~(1)]{KM}, the nonzero module $(\omega_R)_{\p}$ is a canonical module of $R_{\p}$. Since $R_{\p}$ is local, some $R_{\p}$-linear map $(\omega_R)_{\p}\to R_{\p}$ is surjective. Its source and target have the same length, so $R_{\p}$ is Gorenstein.

If $R_{\p}$ is Gorenstein for some prime $\p$ minimal over $I$, then $I^kS_{\p}$ is primary to $\p S_{\p}$. In the regular local ring $S_{\p}$, which contains $\kk$, it has the proper factors $IS_{\p}$ and $I^{k-1}S_{\p}$. By \cite[Theorem~2]{Huneke} its height is at most one. Since $0\ne I\subseteq\p$ and $\operatorname{ht}_S\p=\operatorname{ht}(I^kS_{\p})$, we get $\operatorname{ht}I=1$.

Suppose $\operatorname{ht}I=1$ and set $J=I:_S f$, so $I^k=f^kJ^k$. No height-one prime contains $J$. Multiplication by $f^k$ gives an exact sequence of $S$-modules
\[
0\longrightarrow S/J^k\xrightarrow{\ \cdot f^k\ }R\longrightarrow S/(f^k)\longrightarrow0.
\]
If $J\ne S$, then $J^k$ contains an $S$-regular sequence of length two since it has height at least two. Thus $\Ext_S^i(S/J^k,S)=0$ for $i=0,1$, which also holds if $J=S$. The long exact Ext sequence gives $\omega_R\cong\Ext_S^1(S/(f^k),S)\cong S/(f^k)$ up to shift. An $R$-linear map from this cyclic module is determined by the image of $1$, which may be any element of $0:_R f^k$. Hence
\[
\tr_R(\omega_R)=0:_R f^k=(I^k:_S f^k)/I^k=J^k/I^k.
\]
At a height-one prime $\p$ over $I$ we have $JS_{\p}=S_{\p}$, so $(J^k/I^k)_{\p}=R_{\p}$. Thus the trace is not nilpotent, giving \textup{(i)} $\Rightarrow$ \textup{(iii)}.
\end{proof}

Theorem~\ref{thm:powers} and Proposition~\ref{prop:height-one} give the following criterion for $S/I^k$ to be nearly Gorenstein.

\begin{corollary}\label{cor:trace-powers}
Let $I\subsetneq S$ be a nonzero graded ideal and $k\geq2$, and set $R=S/I^k$. If $\dim R=0$, assume that $k\geq3$ and $\operatorname{char}\kk=0$. The following conditions are equivalent:
\begin{enumerate}
\item\label{cond:trace} $\tr_R(\omega_R)\supseteq\m_R$;
\item\label{cond:nearly} $R$ is Cohen--Macaulay and $\tr_R(\omega_R)\supseteq\m_R$, that is, $R$ is nearly Gorenstein;
\item\label{cond:generated} $I$ is generated by one element.
\end{enumerate}
\end{corollary}

\begin{proof}
Under \textup{(\ref{cond:generated})}, $R$ is a hypersurface. It is Cohen--Macaulay and $\tr_R(\omega_R)=R$, so \textup{(\ref{cond:nearly})} and \textup{(\ref{cond:trace})} hold.

Suppose $I$ is not principal, so $n\geq2$. We have $\m_R/\m_R^2\cong\m/\m^2$ because $I^k\subseteq\m^2$. In height one the ideal $I:_S f$ of Proposition~\ref{prop:height-one} is proper and graded. Thus $I:_S f\subseteq\m$ and the trace is contained in $\m_R^k\subseteq\m_R^2$, ruling out \textup{(\ref{cond:trace})}. For $\operatorname{ht}I\geq2$ and $\dim R>0$, the trace is nilpotent by Proposition~\ref{prop:height-one}. If it contained $\m_R$, then $R$ would have dimension zero. Finally, if $\dim R=0$, the ring is Cohen--Macaulay and $\operatorname{ht}I=n\geq2$. Here Theorem~\ref{thm:powers} gives $\tr_R(\omega_R)\subseteq\m_R^{k-1}\subseteq\m_R^2$, again ruling out \textup{(\ref{cond:trace})}.
\end{proof}

The Artinian case of Question~\ref{question:powers} in positive characteristic leads to the following conjecture.

\begin{conjecture}\label{conj:positive-characteristic}
Suppose that $\operatorname{char}\kk>0$ and $n\geq2$. If $I\subset S$ is a graded $\m$-primary ideal, then $S/I^k$ is not nearly Gorenstein for any $k\geq3$.
\end{conjecture}

\medskip
\noindent\textbf{Acknowledgments.} The author was supported by JSPS KAKENHI Grant Number 25KJ1744.

\medskip
\noindent\textbf{Declaration of AI use.} The author used ChatGPT (OpenAI) to discuss mathematical arguments and examine proofs, search the literature and check references, revise the text, and prepare the \TeX{} source. The author takes full responsibility for the mathematical content, references, and final text.

\end{document}